%% file: garivier.tex
\documentclass[conference]{IEEEtran}
\usepackage{amsmath,amssymb,bbm,ifthen,twoopt,url,amsbsy,epsfig,graphicx}
\usepackage{float,algorithm,algorithmic}
\floatname{algorithm}{Algorithme}

\ifCLASSINFOpdf
\else
\fi
\input{./macro}

\begin{document}
%
\title{Informational Confidence Bounds for Self-Normalized Averages and Applications}

\author{\IEEEauthorblockN{Aur\'elien Garivier}
\IEEEauthorblockA{Institut de Math\'ematiques de Toulouse\\
Universit\'e Paul Sabatier\\
118 route de Narbonne 31062 Toulouse cedex 9\\
Email: aurelien.garivier@math.univ-toulouse.fr}
}


%


\maketitle

\begin{abstract}
We present deviation bounds for self-normalized averages and applications to estimation with a random number of observations.
The results rely on a peeling argument in exponential martingale techniques that represents an alternative to the method of mixture. 
The motivating examples of bandit problems and context tree estimation are detailed.
\end{abstract}


\section{Introduction}
Contrary to a very usual assumption in statistics, some situations require parameter estimation based on samples of \emph{random size}.
Let us first briefly present two probabilistic models of that kind which motivated the derivation of the results presented below.
\subsection{Motivating examples}
\subsubsection{Bandit Problems}
Estimation is sometimes used as a intermediate step in a decision process, and the results can influence the presence of further observations.
Paradigmatic of this situation are \emph{bandit problems}, named in reference to the archetypal
situation of a gambler facing a row of slot-machines and sequentially deciding which one to choose in order to maximize her gains.
The basic model is the following: an agent sequentially chooses actions in a finite set of possible options. Each
action leads to an independent stochastic reward whose distribution is unknown.
What dynamic allocation rule should she choose so as to maximize her
cumulated reward?
Originally motivated by medical trials, this simple model dates back to the
1930s ; it has recently raised a renewed interest because of computer-driven applications, from
computer experiments to recommender systems and Big Data, and numerous variants have been considered (see~\cite{bubeck_cesa-bianchi_2012-ftml} for a recent survey, and \cite{BubeckErnstGarivier11gooducb} for a related model).

One possible solution consists in constructing, at time $t$, a confidence interval based on all past observations for the expected reward associated to each action, then choosing the action with highest upper confidence bounds (UCB).
This rule, popularized by~\cite{AuerEtAl02FiniteTime}, was recently improved and shown to have some optimality properties~\cite{COLT2011,cappeEtal13AOSbandits}.
Obviously, the number of observations used to construct the confidence interval strongly depends on the value of these observations, and standard formulas for fixed-size samples do not apply directly. 
The key element in the recent improvements of this algorithm was the introduction of the informational self-normalized deviation inequalities presented below.

\subsubsection{Context Tree Estimation}
Context tree models, introduced by Jorma Rissanen in \cite{rissanen1983} as efficient tools in Information Theory, have been successfully studied and used since then in many fields of Probability and Statistics, including Bioinformatics, Universal Coding, Mathematical Statistics or Linguistics.
Sometimes also called Variable Length Markov Chain, a context tree process is informally defined as a  Markov chain whose memory length depends on past symbols. 
This property  makes it possible to represent the set of memory sequences as a tree, called the \emph{context tree} of the process. 

A remarkable tradeoff between flexibility and parsimony explains this success: no more difficult to handle than Markov chains, they include memory only where necessary.
Not only do they provide more efficient models for fitting the data: it appears also that, in many applications, the shape of the context tree has a natural and informative interpretation. In Bioinformatics, they have been used to test the relevance of protein families databases \cite{busch2009} and in Linguistics, tree estimation highlights structural discrepancies between Brazilian and European Portuguese~\cite{galves2009}.

Of course, practical use of context tree models requires the possibility of constructing efficient estimators of the model generating the data. Despite the multiplicity of candidate trees, several procedures have been proposed and proved to be consistent, including pruning algorithms~\cite{rissanen1983}, and  \emph{Penalized Maximum Likelihood} estimators (see~\cite{csiszar2006,garivier06bicCons} and references therein, see also~\cite{buhlmann1999}).
These apparently different ideas are in fact closely related~\cite{garivierLeonardi11CTestimation}, the key point being an efficient estimation of the conditional transition probabilities. 
But for a given sample size, the number of transitions observed from a given context is random, and depends on the values of these transitions. Hence, again, sharp deviation bounds for random-sized averages are required in order to obtain efficient memory estimators.

\subsection{Self-Normalized Process}

Several approaches have been proposed to address this problem.
The most obvious is to use a simple union bound on all the possible values of the sample size (as for instance in~\cite{AuerEtAl02FiniteTime}), but this appears to be most often overly pessimistic and significantly sub-optimal.
A more refined treatment consists, when possible, in first lower-bounding the size of the sample, thus upper-bounding the variance of the estimator, and then in using this upper-bound of the variance, for example with Bernstein's maximal inequality for martingales (see e.g.~\cite{celisseEtAl11} for an example on the consistency of a Stochastic Block Model estimator, or \cite{CbLuSt05refinded} on prediction with expert advice).

The most satisfying approach, however, is to consider the associated \emph{self-normalized process}.
From the estimator's point of view, as the size of the sample grows, something changes only when a new observation appears. At those (random) times, the internal clock of the estimator increases by $1$. When the $n$th observation has been reached, the internal clock has a random value which is at most equal to $n$, and on which the variance of the estimator depends. The confidence interval must be constructed accordingly, by taking into account the maximal deviations of the self-normalized deviation process. 

This paper focuses on the case where a sequence of non-asymptotic confidence intervals $[a_t,b_t]$ is required for the common expectation $\mu$ of independent, real-valued random variables $(X_t)_{t\geq 1}$, and where all the confidence intervals are required to be jointly valid over an epoch $t\in\{1,\dots,n\}$. In other words, for all positive $\alpha$, the goal is to construct $\sigma(X_1,\dots,X_t)$-measurable random variables $a_t$ and $b_t$ so as to ensure that the event $\bigcap_{t\leq n} \left\{\mu\in[a_t, b_t]\right\}$ has probability at least $1-\alpha$.

In order to obtain sub-gaussian deviation bounds, the \emph{method of mixture}  (see~\cite{DeLaPena:al:04,pena2009self} and references therein) provides a powerful and elegant tool, recently used in~\cite{AYPS11selfnormBandits} for bandit problems. 
The results presented below follow a different path. Rather than a mixture, they rely on a \emph{peeling} device: the possible numbers of observations are divided into exponentially growing slices, which are treated independently. On each slice, a Cramer-type bound is obtained by a maximal inequality for martingales. 

These results can be considered, in some sense, as non-asymptotic counterparts of the Law of the Iterated Logarithm for martingales. 
They are presented so as to clearly emphasize the cost of the randomness of the sample size: namely, a logarithmic factor of $n$ in front of the exponential Cramer bound (instead of a factor $n$ for the union bound). 
The proof method is generic enough to apply to a large variety of situations, and in particular not only to the sub-gaussian case. 

\subsection{Informational Confidence Bounds}

If the bounds of these confidence intervals are classically 
 chosen to be symmetric around the empirical mean $\Xb_t$, so that $\Xb_t-a_t = b_t-\Xb_t = c/\sqrt{t}$ for a given constant $c$, then the above discussion shows that one needs to control the following supremum of the self-normalized process:
\[ \sup_{t\leq n} \sqrt{t}\left|\Xb_t-\mu \right| \;.\]
This choice, however, is often sub-optimal and was not sufficient in the applications mentioned above. 
The approach used below is somewhat different: the deviations of  $\Xb_t$ are not measured in absolute value, but using a information deviation measure, leading to possibly asymmetric confidence bounds.
Let us recall it briefly: suppose that, for all possible values of the expectation $\mu$, the following Cramer-type inequality with rate function $I(\cdot, \mu)$ is satisfied:
\[\forall x_t\geq \mu, P(\Xb_t \geq x_t)\leq\exp(-t I(x_t; \mu))\;.\]
For a concrete example, one may think about i.i.d. Bernoulli variables with $I(x;\mu) = \kl(x, \mu) = x\log(x/\mu) + (1-x)\log((1-x)/(1-\mu))$.
As the function $I(\cdot; \mu)$ increases on $[\mu, +\infty[$, this bound can be rewritten $P(I(\Xb_t;\mu) \geq I(x_t;\mu), \Xb_t\geq\mu) \leq \exp(-t\, I(x_t;\mu))$
or, defining $\delta = t I(x_t;\mu)$, $P(t\,I(\Xb_t;\mu) \geq \delta, \Xb_t\geq\mu)\leq \exp(-\delta)$;
proceeding similarly on the other side of $\mu$, one obtains 
  \[P\left(t\,I(\Xb_t;\mu)\geq \delta\right) \leq 2 \exp(-\delta)\;.\]
Consequently, one is tempted to choose, as a confidence interval of risk $\alpha$, a neighborhood of $\Xb_t$ in the sense of the pseudo-distance~$I$:
\[[a_t, b_t] = \left\{ \mu : t\,I(\Xb_t ; \mu) \leq \log\frac{2}{\alpha}\right \}\;.\]
Observe that $\mu \in [a_t, b_t]$ if and only if $t\,I(\Xb_t ; \mu) \leq \log\frac{2}{\alpha}$. 
For a \emph{sequential} confidence intervals of this kind, where $P\left(\bigcap_{t\leq n} \left\{\mu\in[a_t, b_t]\right\}\right)$ needs to be controlled, one is thus led to study
 \begin{equation}\label{eq:sni} \sup_{t\leq n} t\, I\!\left( \Xb_t ;\mu \right)\;.\end{equation}

In Section~\ref{sec:sndi},  deviation bounds for \eqref{eq:sni} are presented. The generic result of Theorem~\ref{th:unifbnd:basic} is refined, under some additional hypotheses, in Theorem~\ref{th:unifbnd:logconcave} and Equation~\eqref{eq:unifbnd:subgaussian}. Theorem~\ref{th:unifbnd:basiclli} contains a variant that does not require an upper-bound on the sample size.
A subgaussian inequality is given for the discounted case in Equation~\eqref{eq:unifbnd:discounted}.
In Section~\ref{sec:estimation}, these results are applied to estimation in various models: one-parameter canonical exponential famillies, bounded variables, and multinomial distributions.

\section{Self-Normalized Deviation Inequalities}\label{sec:sndi}
For an increasing filtration $(\mathcal{F}_t)_{t\geq 0}$ on some probability space, consider an adapted, real-valued discrete time process $(S_t)_{t\geq 0}$ such that $S_0=0$. Further assume that the increments $X_t = S_t-S_{t-1}$ are bounded as follows: there exist  $\lambda_1\in[-\infty, 0[$, $\lambda_2\in]0,+\infty]$, and a function $\phi:]\lambda_1,\lambda_2[ \to \R$ such that for all $\lambda\in]\lambda_1,\lambda_2[$ and for all $t\geq 1$:
\[\E\left[ \exp(\lambda X_t) \left| \mathcal{F}_{t-1}\right.\right] \leq \exp\left( \phi(\lambda) \right)\;.\]
In other words, the function $\phi$ dominates the logarithmic moment-generating function (lmgf) of the increments $(X_t)_t$ that are assumed to share the same finite expectation $\mu$.
If the increments $X_t$ are identically distributed, $\phi$ can be chosen as the common lmgf, but it proves useful to consider more general cases. Nevertheless, $\phi$ will be supposed to satisfy all usual properties of a lmgf (see~\cite{DemboZeitouni10largeDevs}, Chapter 2)~: $\phi$ is convex and smooth over $]\lambda_1,\lambda_2[$, $\phi(\mu)=0$; its Legendre transform $I(\cdot;\mu)$,defined on $\R$ as
\[I(x; \mu) = \sup_{\lambda\in\R} \{\lambda x - \phi(\lambda)\}\;,\]
is a convex rate function whose domain is included in $\R^+\cup\{+\infty\}$; it is finite and smooth on an open interval  $\mathcal{D}_I\subset \R$ containing $0$, such that $I(\mu,\mu)=0$. For all $x$ such that $I(x)<\infty$, there exists a unique real number $\lambda(x)\in]\lambda_1,\lambda_2[$ such that
\[\phi'(\lambda(x)) = x \hbox{\quad and }I(x; \mu) = \lambda(x) x-\phi(\lambda(x))  \;.\]
$I(x;\mu)$ tends to infinity with $x$, and can be equal to $+\infty$ outside of some interval $(x_-,x_+)$ where it is finite: it holds that $P(X_t\in[x_-,x_+])=1$, and the limit of $I(\cdot, \mu)$ when tends to $x_+$ is denoted $I_+$.
Under those assumptions, the following result holds:
\begin{theorem}
\label{th:unifbnd:basic}
For every $\delta>0$, 
\begin{multline*}
 P\left(\exists t\in\{1,\dots,n\}: t\,I(\Xb_t;\mu) \geq \delta\right) \\\leq 2e\left\lceil \delta\log(n) \right\rceil \exp(-\delta)\;.\end{multline*}

\end{theorem}
\subsection{Short Proof of Theorem~\ref{th:unifbnd:basic}}
The proof of this result, short enough to be sketched here, is inspired by the proof of the Law of the Iterated Logarithm for martingales that can be found in~\cite{neveu72}. 
The epoch  $\{1,\dots,n\}$ is divided into ``slides'' $\{t_{k-1}+1,\dots,t_k\}$ of exponentially increasing sizes: let $t_0=0$, let $\eta>0$ and, for every positive integer $k$, let $t_k=\left\lfloor(1+\eta)^k\right\rfloor$. 
Denoting $D=\left\lceil\log(n)/\log(1+\eta)\right\rceil$ the smallest integer such that $t_D\geq n$, the union bound yields~:
\begin{equation*}
P\left( \bigcup_{t=1}^n\left\{tI\left( \Xb_t; \mu \right) \geq \delta \right\} \right) \leq \sum_{k=1}^D P(A_k)\;,
\end{equation*}
where $A_k = \bigcup_{t=t_{k-1}+1}^{t_k}\left\{tI\left( \Xb_t; \mu \right) \geq \delta \right\}$.
Denote by $s$ the smallest integer such that $\delta/(s+1)\leq I_+$~: for $t\leq s$, obviously  $P(t\,I(\Xb_t; \mu)\geq \delta, \Xb_t>\mu )  = 0$ and thus $P(A_k)=0$ if $t_k\leq s$. 

Let $k$ be such that $t_k>s$, and $\tti_{k-1} = \max\{t_{k-1}, s\}$.
For every $t\in\{\tti_{k-1}+1,\dots, t_k\}$, there exists $x_t\in[\mu, x_+]$ such that $t\,I(x_t; \mu)=\delta$. 
Let $\lambda_k = \lambda(x_{t_k})$, so that $I(x_{t_k};\mu) = \lambda_kx_{t_k} - \phi(\lambda_k)$, and consider the super-martingale $(W_t^k)_t$ defined by $W_0^k=1$ and, for every $t\geq 1$, $W_t^k = \exp\left( \lambda_k S_t  - t\phi(\lambda_k) \right)\;.$
A maximal inequality ensures that, for all positive real $c$,
\begin{equation*}
 P\left( \bigcup_{t=t_{k-1}+1}^{t_k}\left\{ W_t^k \geq c\right\} \right)
\leq \frac{1}{c}\;.
\end{equation*}
Let us deduce an upper-bound for $P(A_k)$. As $t\,I(x_t; \mu)=\delta$, it holds that
\begin{equation*}
 I(x_{t_k};\mu) \leq I(x_t;\mu) < I(x_{t_k}; \mu)\left( 1+\eta \right)\;.
\end{equation*}
As $I(\cdot; \mu)$ is increasing on the right side of $\mu$, $x_t\geq x_{t_k}$ and
\[\lambda_k x_t - \phi(\lambda_k) \geq \lambda_k x_{t_k} - \phi(\lambda_k) = I(x_{t_k};\mu)\geq \frac{I(x_t; \mu)}{1+\eta}\;.\]
Hence, if $tI\left(\Xb_t; \mu\right)\geq \delta \hbox{ and } \Xb_t\geq \mu$, then $\lambda_k\Xb_t - \phi(\lambda_k) \geq \lambda_k x_t -\phi(\lambda_k)\geq \frac{\delta}{t(1+\eta)}$ and $\lambda_kS_t -t\phi(\lambda_k) \geq \frac{\delta}{1+\eta}$, and thus $W_t^k \geq \exp\left( \frac{\delta}{1+\eta} \right)$. This entails that
\begin{multline*}
 P\left(  \bigcup_{t=t_{k-1}+1}^{t_k} \left\{ tI\left( \Xb_t; \mu \right) \geq \delta\right\}\cap \left\{ \Xb_t > \mu\right\}  \right)\\ \leq
 P\left(  \bigcup_{t=t_{k-1}+1}^{t_k}\left\{W_t^k \geq \exp\left( \frac{\delta}{1+\eta} \right) \right\} \right) 
 \\\leq \exp\left(-\frac{\delta}{1+\eta}  \right)\;.
\end{multline*}
The case $\Xb_t < \mu$ can be treated similarly, and the first claim of the theorem follows. 
The second claim is a consequence of the inequality  $\log(1+1/(\delta-1))\geq 1/\delta$, applied with the approximately optimal choice $\eta=\delta/(\delta-1)$.


Remark that the simple bound of Theorem~\ref{th:unifbnd:basic} highlights the cost of time uniformity: a factor $2e\lceil\delta\log(n)\rceil$, instead as the factor $n$ given by the union bound.
The fact that this cost is sub-polynomial in $n$ appears (especially in~\cite{COLT2011,garivierMoulines11,garivierLeonardi11CTestimation,cappeEtal13AOSbandits}) to be crucial in the analysis of some algorithms and estimators.

\subsection{Improvements and Variants}
This result can be significantly improved under some additional assumptions on the function $I(\cdot; \mu)$:
\begin{theorem}
\label{th:unifbnd:logconcave}
Let $\delta>0$. If the function $I(\cdot; \mu)$ is log-concave, then for every $\eta>0$
\begin{multline*}
 P\left(\exists t\in\{1,\dots,n\}: t\,I(\Xb_t;\mu) \geq \delta\right) \\\leq 2\left\lceil\frac{\log n}{\log \left( 1+\eta \right)}\right\rceil \exp\left(-\left( 1-\frac{\eta^2}{8} \right)\delta\right)\;.
\end{multline*}
In particular, for $\eta = 2/\sqrt{\delta}$, one obtains:
\begin{multline*}P\left(\exists t\in\{1,\dots,n\}: t\,I(\Xb_t;\mu) \geq \delta\right)\\ \leq 2\sqrt{e}\left\lceil \frac{\sqrt{\delta}}{2}\log(n) \right\rceil \exp(-\delta)\;.\end{multline*}
\end{theorem}
The law of the Iterated Logarithm suggests that such a result is hardly improvable: in a sub-gaussian setting where $I(x;\mu)\geq (x-\mu)^2/(2\sigma^2)$, it implies indeed that for all $c>1$:
\begin{multline*}
 P\left(\sup_{t\leq n}\frac{S_t-t\mu}{\sqrt{2\sigma^2t\log\log(n)}}>c\right)\\ \leq 
P\left(\sup_{t\leq n}tI\left( \Xb_t; \mu \right)>c^2\log\log(n)\right) \to 0
\end{multline*}
when $n$ tends to infinity.
Observe that the log-concavity of $I(\cdot,\mu)$, although not always satisfied (even for bounded variables), is reasonable at least locally around $\mu$ if one thinks to the gaussian regime.

Let us mention that in the quadratic (gaussian) case $I(x;\mu) = 2(x-\mu)^2/K^2$, the bound can be slightly improved:
\begin{multline}\label{eq:unifbnd:subgaussian}
P\left(\exists t\in\{1,\dots,n\}: t\,I(\Xb_t;\mu) \geq \delta\right) \\\leq 2\left\lceil\frac{\log n}{\log \left( 1+\eta \right)}\right\rceil \exp\left(-\left( 1-\frac{\eta^2}{16} \right)\delta\right)\;.
\end{multline}
Finally, the method can be adapted in order to obtain non-asymptotic that hold for all $t\geq 1$ in the spirit of the Law of the Iterated Logarithm:
\begin{theorem}
\label{th:unifbnd:basiclli}
For all  $\delta>1$ and all $c>1$,
\begin{multline*}
 P\left(\exists t\geq 1: t\,I(\Xb_t;\mu) \geq \frac{\delta c}{\delta-1}\log\log t + \delta\right) \\\leq \frac{2 \rme c \delta^c}{c-1}\exp(-\delta)\eqsp.
\end{multline*}
In particular, for $c=1+1/\log(\delta)$, one obtains:
\begin{multline*}P\left(\exists t\geq 1: t\,I(\Xb_t;\mu) \geq \frac{\delta(1+\log \delta)}{(\delta-1)\log \delta}\log\log t + \delta\right) \\\leq 2e^2\delta\exp(-\delta)\eqsp.\end{multline*}
\end{theorem}

\subsection{Self-Normalized Form}
In the applications mentioned above, the necessity to guarantee the joint validity of the confidence intervals over the entire epoch comes from the fact that the variables $X_t$ are observed only episodically in a predictable way: there exists, for each  $t\in\{1,\dots,n\}$, a $\{0,1\}$-valued, $\mathcal{F}_{t-1}$-measurable random variable $\varepsilon_t\in\{0,1\}$ such that the current estimate at time $n$ is 
\begin{equation}\label{eq:obsOptionnelle}
\Xb(n) = S(n)/N(n) \,
\end{equation}
where $S(n) =\sum_{t=1}^n \varepsilon_t X_t$ and $N(n) = \sum_{t=1}^n \varepsilon_t$.
Theorem~\ref{th:unifbnd:basic} yields:
\[P\left(I\left(\Xb(n); \mu \right) \geq \frac{\delta}{N(n)}  \right)\leq 2e\left\lceil\delta\log(n)\right\rceil \exp(-\delta)\;.\]

In the definition~\eqref{eq:obsOptionnelle}, $S(n)$ is written as a martingale transform or, equivalently, a discrete stochastic integral. Continuous-time variants of Theorem~\ref{th:unifbnd:basic} can be obtained following the same lines (using the same peeling trick) for stochastic integrals.

Furthermore, this approach can be adapted to non-stationary contexts: assume for simplicity 
that the variables $(X_t)_t$ are independent, of expectation $\mu_t$ respectively, and that their absolute value is almost-surely bounded by $B$. If $\mu_t$ does not change too fast (or too often) with $t$, one may consider the discounted estimator $\Xb_\gamma(n)$ of  $\mu_n$ defined by
\[\Xb_\gamma(n) = \frac{S_\gamma(n)}{N_\gamma(n)}\;,\]
where  $\gamma\in]0,1[$,  $S_\gamma(n) =\sum_{t=1}^n \gamma^{n-t}\varepsilon_t X_t$ and $N_\gamma(n) = \sum_{t=1}^n \gamma^{n-t}\varepsilon_t$.
The difference between  $\Xb_\gamma(n)$ and $\mu_n$ can be decomposed into a term of bias (which is not discussed here) and a fluctuation term $\Xb_\gamma(n) - M_\gamma(n)/N_\gamma(n)$, where $M_\gamma(n) = \sum_{t=1}^n\gamma^{n-t}\varepsilon_t\mu_t$.
This fluctuation term can be controlled by the adapting the martingale techniques above: one obtains that
\begin{multline}\label{eq:unifbnd:discounted}
 P\left( \frac{S_\gamma(n)-M_\gamma(n)}{\sqrt{N_{\gamma^2}(n)}} \geq \delta \right) \\\leq \left\lceil\frac{\log \nu_\gamma(n)}{\log(1+\eta)}\right\rceil \exp\left(-\frac{2\delta^2}{B^2}\left( 1-\frac{\eta^2}{16} \right) \right)\;,
\end{multline}
with $\nu_\gamma(n) = \sum_{t=1}^n \gamma^{n-t} = (1-\gamma^n)/(1-\gamma)<\min\{(1-\gamma)^{-1}, n\}$.
This results permits to analyze the \emph{Discounted-UCB} algorithm~\cite{garivierMoulines11} earlier proposed by Kocsis Szepesv\'ari~\cite{kocsis:szepesvari:2006}.

\section{Application to Estimation}\label{sec:estimation}
Let us now show briefly how these inequalities may be used in the analysis of some stochastic algorithms. 
The key point is that Theorem~\ref{th:unifbnd:basic} allows the construction of a sequence of confidence intervals  $([a_t,b_t])_{1\leq t\leq n}$ for $\mu$ that are simultaneously valid with high probability. The interval
\[[a_t, b_t] = \left\{\mu \in\left[x_-,x_+\right]: tI\left(\Xb_t; \mu\right) \leq \delta\right\} \]
contains all the values in a neighborhood of  $\Xb_t$ in the sense of the pseudo-distance defined by $I$. By Theorem~\ref{th:unifbnd:basic},
\[P\left( \bigcap_{t=1}^n \left\{ \mu\in [a_t,b_t]\right\} \right) \geq 1- 2e\left\lceil \delta\log(n) \right\rceil \exp(-\delta)\;.\]
Similarly, one obtains obtains confidence intervals for the case presented in Equation~\eqref{eq:obsOptionnelle}. 
This framework applies as well in bandit problems, where only the reward of the chosen arm is observed, that the estimation of Markovian models where, at each time, only the estimates relative to the current past observations are updated.
Of course, in these examples, the identity of the variable(s) observed at time $t$ is absolutely not independent of the past observations.
By choosing $\delta$ such that  $2e\left\lceil \delta\log(n) \right\rceil \exp(-\delta)\leq \alpha$, one obtains the confidence interval  $\left\{\mu : I\left(\Xb(n); \mu\right) \leq \delta/N(n)\right\} $ of risk at most $\alpha$.

\subsection{One-Parameter Exponential Model and Bounded Variables}
In this section, we assume that the variables $(X_t)_t$ are independent and identically distributed, and that their distribution $P_{\theta_0}$ belongs to a canonical exponential model of the form $\{P_\theta : \theta\in\Theta\}$, where $\Theta$ is an real interval and where $P_\theta$ has, with respect to some reference measure, the density $p_\theta:\R\to\R$ defined by:
\[p_\theta(x) = \exp\left( x\theta - b(\theta) + c(x) \right)\;.\]
Here, $c$ is a real function and the log-partition function $b$ is supposed to be twice differentiable.
It is well-known that, by denoting  $\mu(\theta)=\dot{b}(\theta)$ the expectation of $P_\theta$, one defines a one-to-one, differentiable mapping $\mu$.
In this case, one easily shows that the rate function  $I$ is directly related to the  Kullback-Leibler divergence (which is here a Bregman divergence for $b$) as follows: for every $\beta,\theta\in\Theta$, 
 \[\KL(P_\beta; P_\theta) = I(\mu(\beta); \mu(\theta)) = b(\theta)-b(\beta) - \dot{b}(\beta) (\theta-\beta) \;.\]
Hence, a sequence $(R_t)_{t\geq 1}$ of confidence intervals for the parameter $\theta_0$ jointly valid with probability  $1- 2e\left\lceil \delta\log(n) \right\rceil \exp(-\delta)$ is obtained by choosing:
\begin{align*}
R_t &= \left\{ \theta : \KL\left(P_{\mu^{-1}(\Xb_t)}; P_{\theta}  \right) \leq \frac{\delta}{t}\right\} \\
&= 
\left\{ \theta : I\left(\Xb_t; \mu(\theta)  \right) \leq \frac{\delta}{t}\right\}
\;.\end{align*}
This applies in particular to usual families of distributions like Poisson, Exponential, Gamma (with fixed shape parameter)...
In~\cite{COLT2011}, an example concerning exponential variable detailed: in that case, $I(x,y) = x/y-1-\log(x/y)$.

But the case of Bernoulli variables deserves to be highlighted, as it easily extends to general \emph{bounded} variables.
Indeed, as observed by Hoeffding~\cite{Hoeffding:63}, the exponential moments of a $[0,1]$-valued variable $X$ with expectation $\mu$ are upper-bounded by those of a Bernoulli variable, and for all $\lambda \in \R$ it holds that
\[E\left[ \exp(\lambda X) \right] \leq 1-\mu +\mu\exp(\lambda) \;,\]
with equality if and only if $X\sim\mathcal{B}(\mu)$.
Recall that $\kl$ denotes the binary entropy function, i.e. the rate function associated to Bernoulli variables. Theorem~\ref{th:unifbnd:basic} yields that, for independent variables  $X_t$ bounded in $[0,1]$, 
\begin{equation}\label{eq:ThBounded}P\left( \sup_{t\leq n} \kl\left( \Xb_t, \mu\right)\geq \frac{\delta}{t} \right) \leq 2e\left\lceil \delta\log(n) \right\rceil \exp(-\delta)\;.\end{equation}
Of course, this result together with Pinsker's inequality $\kl(p,q)\geq 2(p-q)^2$, yields a self-normalized version of Hoeffding's inequality on the epoch $t\in\{1,\dots,n\}$: 
\begin{equation}\label{eq:ThBoundedHoeffding}
 P\left( \sup_{t\leq n} \left| \Xb_t - \mu\right|\geq \frac{\delta}{\sqrt{t}} \right) \leq 4e\left\lceil \delta^2\log(n) \right\rceil \exp(-2\delta^2)\;.
\end{equation}
This bound may seem simpler and easier to use than the previous one. 
However, the case of bounded bandits (detailed in~\cite{COLT2011,cappeEtal13AOSbandits}), as well as the case of context tree estimation (presented in~\cite{garivierLeonardi11CTestimation}) show that Equation~\eqref{eq:ThBounded} is sometimes really to be preferred, as it leads to significantly more efficient algorithms at the price of an hardly increased computational complexity.

\subsection{Multinomial Distributions}
As suggested by Sanov's (asymptotic) Theorem, this kind of inequalities is not limited to real-valued variables. 
It is also possible to construct informational, self-normalized confidence regions for random vectors; let us detail here the simple case of multinomial laws, as they are required, for example, in order to estimate transition distributions in Markov chains (see~\cite{garivierLeonardi11CTestimation, allerton10}).
Let $P$ and $Q$ be two elements of the set  $\mathcal{S}$ of all probability distributions over a finite set $A$. 
By remarking that 
\begin{align*}
-&\KL(P;Q) +\sum_{x\in A}\kl\left(P(x);Q(x)\right)  \\
 &=  \left(|A|-1\right) \sum_{x\in A}\frac{1-P(x)}{|A|-1}\log\left( \frac{(1-P(x))/(|A|-1)}{(1-Q(x))/(|A|-1)} \right)
\end{align*}
is non-negative, one easily shows that 
 \begin{equation*}
\KL(P; Q) \leq \sum_{x\in A}\kl\left(P(x);Q(x)\right)\,.
\end{equation*}
It follows that if $X_1,\dots,X_n$ are i.i.d. variables of law $P_0\in\mathcal{S}$, and if  $\Ph_t(k) = \sum_{s=1}^t \1\{X_s=k\}/t$, then
\begin{align}
P\Big(&\exists t\in\{1,\dots,n\} : \;\;\KL\left( \Ph_t; P_0 \right) \geq \frac{\delta}{t} \Big) \nonumber\\
&\leq  \sum_{a\in A} P\left(\exists t\in\{1,\dots,n\}  : \kl\left( \Ph_t(a); P_0(a) \right)\geq\frac{\delta}{|A|t} \right)\nonumber\\
&\leq 2e\left( \delta \log(n)+|A| \right) \exp\left( -\frac{\delta}{|A|} \right)\;.\label{eq:borneKL}
\end{align}

The fact that this bound involves directly the Kullback-Leibler divergence between the empirical measure and the true distribution allows, in context tree estimation (see~\cite{garivierLeonardi11CTestimation}), to suppress unnecessary assumptions that resulted, in previous papers, from the use of Bernstein's inequality.
Moreover, the Equation~\eqref{eq:borneKL} permits to construct a sequence $(R_t)_{t\leq n}$ of ``Sanov-type'' confidence regions for $P_0$  that are simultaneously valid with probability at least $1-\alpha$, by choosing Kullback-Leibler neighborhoods of the maximum likelihood estimator:
\[R_t = \left\{ Q \in \mathcal{S} : \KL(\Ph_t; Q) \leq \frac{\delta}{t}\right\}\;,\]
with $\delta $ such that $2e\left( \delta \log(n)+|A| \right) \exp\left( -\delta/|A| \right) =\alpha$.
These regions $R_t$ of the simplex have nice geometric properties that are exploited in~\cite{allerton10} for reinforcement learning in Markov Decision Process, improving on former results using $L^1$ regions.


%
\IEEEpeerreviewmaketitle

\bibliographystyle{IEEEtran}
\bibliography{IEEEabrv,garivier.bib}
%



\end{document}

%% file: macro.tex
 \newtheorem{theorem}{Theorem}

\newcommand{\R}{{\mathbb{R}}}
\newcommand{\E}{{\mathbb{E}}}

\def\lg2{\log_2}
\def\R{\mathbb{R}}     
\def\1{\mathbbm{1}}   \let\geq\geqslant   \let\leq\leqslant
\def\E{\mathbb{E}} 
     
\newcommand{\deftitle}[1]{\begin{definition}{\bf (#1)}}

\newcommand{\assurepasmath}[1]{\ifmmode\text{#1}\else #1\fi}

\newcommand{\kamb}[1]{\ifmmode\text{AFA}_{#1}\else AFA\ensuremath{_{\text{#1}}}\fi}

\newcommand{\ckamb}[1]{\ifmmode\text{{AFC}}_{#1}\else {AFC}\ensuremath{_{\text{#1}}}\fi}

\newcommand{\Sa}{\ifmmode\mathcal{S}\else\ensuremath{\mathcal{S}}\fi}

\newcommand{\auset}[2]{\ifmmode\boldsymbol{D}_{#1}(#2)\else\ensuremath{\boldsymbol{D}_{\text{#1}}(\text{#2})}\fi}
\newcommand{\partitions}[1]{\ifmmode\boldsymbol{P}(#1)\else\ensuremath{\boldsymbol{P}}(#1)\fi}
\newcommand{\partitionsr}[2]{\ifmmode\boldsymbol{P}_{#1}(#2)\else\ensuremath{\boldsymbol{P}_{#1}}(#2)\fi}

\newcommand{\Acc}[2]{\ifmmode Acc_{#1}(#2)\else\ensuremath{Acc_{\text{#1}}(\text{#2})}\fi}

\newlength\jataille
\newcommand{\figgauche}[3]%
{\jataille=\textwidth\advance\jataille by -#1
\parbox{#1}{#2}
\parbox{\jataille}{#3}
}

\newcommand{\figtxt}[1]%
{{\it {\small #1}}}

\makeatletter
\def\hlinewd#1{%
\noalign{\ifnum0=`}\fi\hrule \@height #1 %
\futurelet\reserved@a\@xhline}
\makeatother

\def\rme{\operatorname{e}}

\newcommand{\MReward}[3][{}]{\ifthenelse{\equal{#1}{}}{\bar{X}_{#2}(#3)}{\bar{X}_{#2}(#3, #1)}}
\newcommandtwoopt{\Nplay}[4][{}][{}]{
  \ifthenelse{\equal{#1}{}}{
    \ifthenelse{\equal{#2}{}}
    {N_{#3}(#4)}
    {N_{#3}(#4, #2)}
  }{
    \ifthenelse{\equal{#1}{sum}}{
      \ifthenelse{\equal{#2}{}}
      {n_{#3}(#4)}
      {n_{#3}(#4, #2)}
    }{
      \ifthenelse{\equal{#2}{}}
        {N^{#1}_{#3}(#4)}
        {N^{#1}_{#3}(#4, #2)}
    }
  }
}

\newcommandtwoopt{\Np}[3][{}][{}]{
\ifthenelse{\equal{#1}{}}{
  \ifthenelse{\equal{#2}{}}
    {N(#3)}
    {N^{#2}(#3)}
}{
  \ifthenelse{\equal{#2}{}}
    {N_{#1}(#3)}
    {N_{#1}^{#2}(#3)}}
}
\def\eqsp{\;}

\def\Xb{\bar{X}}
\def\Ph{\hat{P}}
\def\tti{\tilde{t}}

\DeclareMathOperator*{\KL}{\operatorname{KL}}
\DeclareMathOperator*{\kl}{\operatorname{kl}}
